\documentclass[10pt]{article}
\usepackage{lmodern}
\usepackage{amsmath}
\usepackage[T1]{fontenc}
\usepackage[utf8]{inputenc}
\usepackage{authblk}
\usepackage{amsfonts}
\usepackage{graphicx}
\usepackage{tkz-euclide}
\tikzset{elegant/.style={smooth,thick,samples=50,cyan}}
\usepackage{color,tikz}
\usetikzlibrary{decorations.pathreplacing,calc}
\usepackage{rotating}
\usepackage{amssymb}
\usepackage[english]{babel}
\usepackage{color}
\usepackage{amsthm}
\usepackage{graphicx}
\usepackage{mathrsfs}
\usepackage{makecell}
\usepackage{microtype}
\usepackage{enumerate}
\usepackage{mathscinet}
\usepackage{array}
\usepackage{multirow}
\usepackage{booktabs}
\usepackage{enumerate}
\usepackage[cal=boondoxo,bb=ams]{mathalfa}
\usepackage{hyperref}
\hypersetup{hidelinks}
\usepackage{titlesec}
\newtheorem{theorem}{Theorem}[section]
\newtheorem{prop}{Proposition}[section]

\newtheorem{remark}{Remark}[section]

\newcommand{\ml}{\mathcal}
\newcommand{\mb}{\mathbb}

\DeclareMathOperator{\non}{non}
\DeclareMathOperator{\lin}{lin}
\DeclareMathOperator{\intt}{int}

\title{On the critical exponent and sharp lifespan estimates for semilinear damped wave equations with data from Sobolev spaces of negative order} 
\author{Wenhui Chen\thanks{Wenhui Chen (wenhui.chen.math@gmail.com)}}
\author[2]{Michael Reissig\thanks{Michael Reissig (reissig@math.tu-freiberg.de)}}
\affil[1]{School of Mathematical Sciences, Shanghai Jiao Tong University, 200240 Shanghai, China}
\affil[2]{Institute of Applied Analysis, Faculty of Mathematics and Computer Science, Technical University Bergakademie Freiberg, 09596 Freiberg, Germany}

\date{}
\begin{document}
\maketitle
\begin{abstract}
We study semilinear damped wave equations with power nonlinearity $|u|^p$ and initial data belonging to Sobolev spaces of negative order $\dot{H}^{-\gamma}$. In the present paper, we obtain a new critical exponent $p=p_{\mathrm{crit}}(n,\gamma):=1+\frac{4}{n+2\gamma}$ for some $\gamma\in(0,\frac{n}{2})$ and low dimensions in the framework of Soblev spaces of negative order. Precisely, global (in time) existence of small data Sobolev solutions of lower regularity is proved for $p>p_{\mathrm{crit}}(n,\gamma)$, and blow-up of weak solutions in finite time even for small data if $1<p<p_{\mathrm{crit}}(n,\gamma)$. Furthermore, in order to more accurately describe the blow-up time, we investigate sharp upper bound and lower bound estimates for the lifespan in the subcritical case.

\noindent\textbf{AMS Classification (2020)} Primary: 35L71; Secondary: 35L15, 35B33, 35A01, 35B44.
\end{abstract}
\fontsize{12}{15}
\selectfont

\section{Introduction}
In this paper, we investigate the following Cauchy problem for semilinear wave equation with classical damping and power nonlinearity:
\begin{align}\label{Semilinear_Damped_Waves}
	\begin{cases}
		u_{tt}-\Delta u+u_t=|u|^p,&x\in\mb{R}^n,\ t>0,\\
		u(0,x)=\epsilon u_0(x),\ u_t(0,x)=\epsilon u_1(x),&x\in\mb{R}^n,
	\end{cases}
\end{align}
with $p>1$, where the initial data with its size parameter $\epsilon>0$ belongs additionally to some Sobolev spaces of negative order. Our main purpose is to determine a critical exponent for the nonlinear Cauchy problem \eqref{Semilinear_Damped_Waves} with $(u_0,u_1)\in \dot{H}^{-\gamma}\times \dot{H}^{-\gamma}$ equipping $\gamma>0$, in which the Sobolev spaces of negative order are defined by (see, for instance, \cite{Runst-Si-book})
\begin{align}\label{Negative_Sobolev}
	\dot{H}^{-\gamma}:=\Big\{f:\ f\in\ml{Z}'\ \ \mbox{such that}\ \  \|f\|_{\dot{H}^{-\gamma}}:=\|(-\Delta)^{-\frac{\gamma}{2}}f\|_{L^2}=\|\ml{F}^{-1}(|\xi|^{-\gamma}\hat{f})\|_{L^2}<\infty\Big\}.
\end{align}
Here, $\ml{Z}'$ stands for the topological dual space to the subspace of the Schwartz space $\ml{S}$ consisting of function with $\mathrm{d}_{\xi}^k\hat{f}(0)=0$ for all $k\in\mb{N}$, namely, $\ml{Z}'$ is the factor space $\ml{S}'/\ml{P}$, where $\ml{P}$ is the space of all polynomials. The critical exponent means the threshold condition on the exponent $p$ for global (in time) Sobolev solutions and blow-up of local (in time) weak solutions with small data. To be specific, under additional $\dot{H}^{-\gamma}$ assumptions for the initial data, the new critical exponent, which will be proposed for \eqref{Semilinear_Damped_Waves}, is
\begin{align}\label{Critical_Exp}
	p=p_{\mathrm{crit}}(n,\gamma):=1+\frac{4}{n+2\gamma}\ \ \mbox{with}\ \ \gamma\in\Big(0,\min\Big\{\frac{n}{2},\frac{\sqrt{n^2+16n}}{4}-\frac{n}{4}\Big\}\Big),
\end{align}
for $n=1,\dots,6$. A further purpose of this paper is to derive sharp lifespan estimates for weak solutions to the semilinear Cauchy problem \eqref{Semilinear_Damped_Waves}, in which the lifespan $T_{\epsilon}$ of solution is defined by
\begin{align*}
	T_{\epsilon}:=\sup\big\{T>0:\ &\mbox{there exists a unique local (in time) solution }u\mbox{ to the Cauchy problem \eqref{Semilinear_Damped_Waves}}\\
	& \mbox{on }[0,T)\mbox{ with a fixed parameter }\epsilon>0\big\}.
\end{align*}
Taking initial data from $\dot{H}^{-\gamma}$ with $\gamma\in(0,\min\{\frac{n}{2},2\})$ and $1+\frac{2\gamma}{n}\leqslant p\leqslant\frac{n}{(n-2)_+}$, we will demonstrate the sharpness of new lifespan estimates
\begin{align*}
	T_{\epsilon}\begin{cases}
		=\infty&\mbox{if}\ \ p>p_{\mathrm{crit}}(n,\gamma),\\
		\simeq C\epsilon^{-\frac{2}{2p'-2-\frac{n}{2}-\gamma}}&\mbox{if}\ \ p<p_{\mathrm{crit}}(n,\gamma),
	\end{cases}
\end{align*}
where $C$ is an independent of $\epsilon$, positive constant. In the last exponent, $p'$ is the conjugate of $p$ such that $\frac{1}{p}+\frac{1}{p'}=1$.
\color{black}

For the classical semilinear damped wave equation \eqref{Semilinear_Damped_Waves} with initial data belonging additionally to $L^1$ space, it has been proved in \cite{Matsumura-1976,Todorova-Yordanov-2001,Zhang-2001,Ikehata-Tanizawa-2005} that the critical exponent is the so-called \emph{Fujita exponent}
\begin{align*}
	p_{\mathrm{Fuj}}(n):=1+\frac{2}{n},	
\end{align*}
which is also the critical exponent for the semilinear heat equation (see, for example, \cite{Fujita-1966}) as follows:
\begin{align*}
	\begin{cases}
		w_{t}-\Delta w=|w|^p,&x\in\mb{R}^n,\ t>0,\\
		w(0,x)=w_0(x),&x\in\mb{R}^n.
	\end{cases}
\end{align*}
This effect is motivated by the diffusion phenomenon (see, for example, \cite{Marcati-Nishihara-2003,Narazaki-2004}) between the linear classical damped wave equation and the heat equation. Let us return our discussions to the critical exponent for the classical damped wave equation \eqref{Semilinear_Damped_Waves}. In the pioneering paper \cite{Matsumura-1976}, the author proved global (in time) existence of small data solutions for the supercritical case $p>p_{\mathrm{Fuj}}(n)$ if $n=1,2$. Later, the authors of \cite{Todorova-Yordanov-2001} demonstrated global existence for any $n\geqslant 1$ in the supercritical case  by assuming compactly supported initial data, and blow-up of local (in time) solutions in the subcritical case $1<p<p_{\mathrm{Fuj}}(n)$ under a suitable sign assumption for initial data. The blow-up result for the critical case $p=p_{\mathrm{Fuj}}(n)$ was obtained in \cite{Zhang-2001} by using the test function method. Finally, concerning the supercritical case, the restriction of compactly supported data was removed by \cite{Ikehata-Tanizawa-2005}. Quite recently, the authors of \cite{Ebert-Girardi-Reissig-2020} considered the equation \eqref{Semilinear_Damped_Waves} with nonlinear term $\omega (|u|)|u|^{p_{\mathrm{Fuj}}(n)}$ instead of $|u|^p$, where $\omega=\omega (s)$ denotes a suitable modulus of continuity, in which they found sharp conditions for the critical regularity of the nonlinear term. Concerning lifespan estimates for the Cauchy problem \eqref{Semilinear_Damped_Waves} with additional $L^1$ data, the authors of  \cite{LS-Li-Zhou=1995,LS-Kirane-Qafsaoui=2002,LS-Nishihara=2011,LS-Ikeda-Wakasugi=2015,LS-Ikeda-Ogawa=2016,LS-Fujiwara-Ikeda-Wakasugi=2019,LS-Ikeda-Sobajima=2019,LS-Lai-Zhou=2019} derived sharp lifespan estimates as follows:
\begin{align}\label{Life_span}
	T_{\epsilon}\begin{cases}
		=\infty &\mbox{if}\ \ p>p_{\mathrm{Fuj}}(n),\\
		\simeq \exp\big(C\epsilon^{-(p-1)}\big)&\mbox{if}\ \ p=p_{\mathrm{Fuj}}(n),\\
		\simeq C\epsilon^{-\frac{2(p-1)}{2-n(p-1)}}&\mbox{if} \ \ p<p_{\mathrm{Fuj}}(n),
	\end{cases}
\end{align}
where $C$ is a positive constant independent of $\epsilon$.
\color{black}
In recent years, critical exponents for the classical semilinear damped wave equation in other frameworks catch a lot of attentions. For the Cauchy problem \eqref{Semilinear_Damped_Waves} with initial data belonging additionally to $L^m$ spaces and $m\in(1,2)$, the critical exponent is changed into the \emph{modified Fujita exponent} $p_{\mathrm{Fuj}}(\frac{n}{m})=1+\frac{2m}{n}$. Particularly, we stress that differently from the case of $L^1$ data, the global (in time) solution with additional $L^m$ regular data exists uniquely in the critical case $p_{\mathrm{Fuj}}(\frac{n}{m})$. The studies and verification of $p_{\mathrm{Fuj}}(\frac{n}{m})$ to be the critical exponent are shown in \cite{Nakao-Ono,Ikehata-Oho-2002,Narazaki-Nishihara-2008,Ikeda-Inui-Wakasugi-2017,Ikeda-Inui-Okamoto-Wakasugi-2019} and references therein. Concerning semilinear damped wave equation in non-Euclidean frameworks,  the authors of \cite{Georgiev-Palmieri-2020} found the critical exponent $p_{\mathrm{Fuj}}(\ml{Q})=1+\frac{2}{\ml{Q}}$ in the Heisenberg group framework $\mathrm{H}_n$ with the homogeneous dimensions $\ml{Q}:=2n+2$ of $\mathrm{H}_n$ (basing on decay estimates for the linearized problem \cite{Palmieri-JFA-2020}), and the author of \cite{Palmieri-2021} obtained the critical exponent $p_{\mathrm{Fuj}}(0):=\infty$ in the compacted Lie group framework $\mb{G}$. Nevertheless, to the best of authors' knowledge, the framework for initial data localizing in Sobolev spaces of negative order to damped wave equations has never been considered in the literature, even for the linear Cauchy problem. For these reasons, it seems interesting to study qualitative properties (well-posedness, blow-up criterion and behaviors, decay rate, asymptotic profiles) of solutions to  damped wave equations with initial data taken additionally from $\dot{H}^{-\gamma}$ carrying $\gamma>0$.

The main concerns of this paper are to verify the critical exponent and to derive sharp lifespan estimates for weak solutions to the Cauchy problem \eqref{Semilinear_Damped_Waves} with initial data taken additionally from $\dot{H}^{-\gamma}$, where the main results and their explanations are presented in Section \ref{Section_Main_Result}.
Motivated by the papers \cite{Guo-Wang-2012,Guo-Tice=2013}, in Section \ref{Section_Linear}  we establish decay estimates and a diffusion phenomenon for the linear classical damped wave equation with vanishing right-hand side by applying WKB analysis and Fourier analysis. Basing on these results, we construct a time-weighted Sobolev space of lower order for Sobolev solutions to the semilinear Cauchy problem \eqref{Semilinear_Damped_Waves}. Then, we demonstrate global (in time) well-posedness for this semilinear Cauchy problem if $p>p_{\mathrm{crit}}(n,\gamma)$ in Section \ref{Section_GESDS} with the help of Banach's fixed point argument, where the treatment of the nonlinear term in $\dot{H}^{-\gamma}$ is based on some applications of Hardy-Littlewood-Sobolev inequality and the fractional Gagliardo-Nirenberg inequality. Later, in Section \ref{Section_Blow-up}, to conclude the optimality, i.e. the critical exponent \eqref{Critical_Exp}, by employing the test function method, we obtain blow-up of weak solutions even for small data provided $1<p<p_{\mathrm{crit}}(n,\gamma)$, and upper bound estimates for the lifespan as byproduct of our approach. Finally, with the aid of a \emph{contrario}, we derive lower bound estimates for the lifespan of mild solutions in Section \ref{Section_Lower_Bound} to assert the sharpness of our lifespan estimates.

\color{black}

\medskip

\noindent\textbf{Notation:}  Throughout this paper, $c$ and $C$ denote some positive constants, which may be changed
from line to line. We denote that $f\lesssim g$ if there exists a positive
constant $C$ such that $f\leqslant Cg$. Moreover, $f\simeq g$ means that $f\lesssim g$ and $g\lesssim f$. The Japanese bracket is denoted by $\langle x\rangle:=\sqrt{1+|x|^2}$. The Sobolev spaces of negative order $\dot{H}^{-\gamma}$ with $\gamma>0$ were defined in \eqref{Negative_Sobolev}. Finally, we define $(x)_+:=\max\{x,0\}$ and $\frac{1}{(x)_+}=\infty$ when $x\leqslant 0$.\color{black}

\section{Main results and their explanations}\label{Section_Main_Result}
\subsection{Results and discussions for the critical exponent}
Let us state the global (in time) existence result firstly.
\begin{theorem}\label{Thm_GESDS_Lower}
	Let $s\in(0,1]$ and $\gamma\in (0,\frac{n}{2})$ for $n\geqslant 1$. Let the exponent $p$ fulfil
	\begin{align}\label{Condition_p}
		p\begin{cases}
			>p_{\mathrm{crit}}(n,\gamma)&\mbox{if}\ \ \gamma\leqslant \tilde{\gamma},\\
			\geqslant \displaystyle{1+\frac{2\gamma}{n}}&\mbox{if}\ \ \gamma> \tilde{\gamma},
		\end{cases}
	\end{align}
	and  $1<p\leqslant \frac{n}{n-2s}$ if $n>2s$, where $\tilde{\gamma}$ denotes the positive root of $2\tilde{\gamma}^2+n\tilde{\gamma}-2n=0$. Let us assume
	\begin{align*}
		(u_0,u_1)\in \ml{A}_s:=(H^s\cap \dot{H}^{-\gamma})\times(L^2\cap \dot{H}^{-\gamma}).
	\end{align*}
	Then, there exists a constant $C>0$ and a small parameter $0<\epsilon\ll 1$ such that for all $\|(u_0,u_1)\|_{\ml{A}_s}\leqslant C$, \color{black}there is a uniquely determined Sobolev solution
	\begin{align*}
		u\in\ml{C}([0,\infty),H^{s})
	\end{align*}
	to the Cauchy problem for the semilinear damped wave equation \eqref{Semilinear_Damped_Waves}. Hence, the lifespan of solution is given by $T_{\epsilon}=\infty$. Furthermore, the solution satisfies the following estimates:
	\begin{align*}
		\|u(t,\cdot)\|_{L^2}&\lesssim \epsilon (1+t)^{-\frac{\gamma}{2}}\|(u_0,u_1)\|_{\ml{A}_s},\\
		\|u(t,\cdot)\|_{\dot{H}^s}&\lesssim \epsilon (1+t)^{-\frac{s+\gamma}{2}}\|(u_0,u_1)\|_{\ml{A}_s}.
	\end{align*}
\end{theorem}
\begin{remark}
	As we will see later, the derived decay rates in Theorem \ref{Thm_GESDS_Lower} coincide with those of the corresponding linearized damped wave equation \eqref{Damped_Waves} in Proposition \ref{Prop_Decay_Linear}, which verifies the effect of no loss of decay.
\end{remark}
\begin{remark}
	Let us give some examples for the admissible range of exponents $p$ for the global (in time) existence result in low dimensions $n=1,\dots,4$.
	\begin{itemize}
		\item When $n=1,2$, we take $\gamma\in(0,\frac{n}{2})$, $s\in(0,1]$ and the exponent satisfies
		\begin{align*}
			1+\frac{4}{n+2\gamma}<p\begin{cases}
				<\infty&\mbox{if}\ \ n\leqslant 2s,\\[0.5em]
				\displaystyle{\leqslant\frac{n}{n-2s}}&\mbox{if}\ \ n> 2s.
			\end{cases}
		\end{align*}
		\item When $n=3,4$, we take $s\in(0,1]$ and the exponent satisfies
		\begin{align*}
			1+\frac{4}{n+2\gamma}<p\leqslant\frac{n}{n-2s}&\ \ \mbox{if}\ \ 0<\gamma\leqslant\tilde{\gamma},\\
			1+\frac{2\gamma}{n}\leqslant p\leqslant\frac{n}{n-2s}&\ \ \mbox{if}\ \ \tilde{\gamma}<\gamma<\frac{n}{2}.
		\end{align*}
	\end{itemize}
	Moreover, we should emphasize that the positive root $\tilde{\gamma}<2$ for all $n\geqslant 1$.
\end{remark}

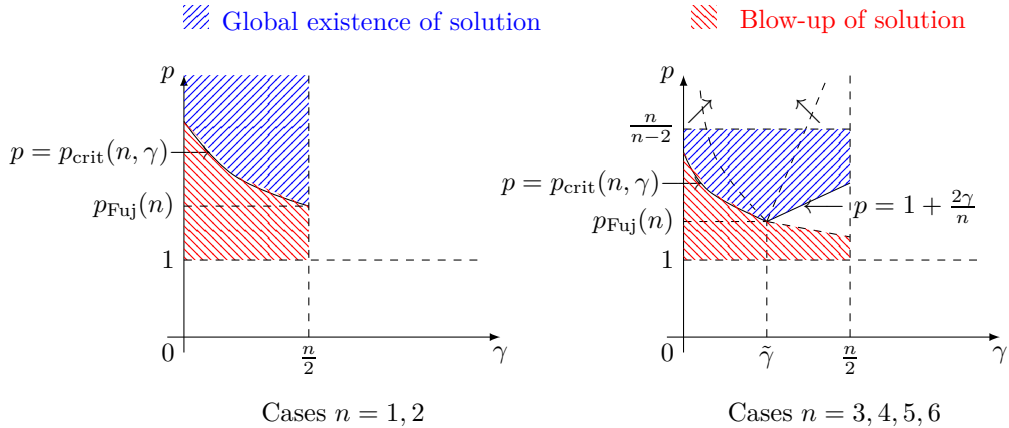
\begin{figure}[http]
	\centering
	\begin{tikzpicture}[>=latex,xscale=1.35,yscale=1.25,scale=0.82]
		\node[left] at (4.4,4.1) {{\color{blue} Global existence of solution}};
		\node[left] at (9.5,4.1) {{\color{red} Blow-up of solution}};
		\fill [pattern=north east lines, pattern color=blue]  (0,4) -- (0, 4.3) -- (0.3,4.3) -- (0.3, 4);
		\fill [pattern=north west lines, pattern color=red]  (6.1,4) -- (6.1, 4.3) -- (6.4,4.3) -- (6.4, 4);
		\draw[->] (-0.2,0) -- (3.8,0) node[below] {$\gamma$};
		\draw[->] (0,-0.2) -- (0,3.4) node[left] {$p$};
		\node[right, color=black] at (-2.3,2.4) {{ $p=p_{\mathrm{crit}}(n,\gamma)$$\longrightarrow$}};
		\node[left] at (0,-0.2) {{$0$}};
		\draw[color=black] plot[smooth, tension=.7] coordinates {(0,2.8) (0.6,2.1) (1.5,1.7)};
		\fill [pattern=north west lines, pattern color=red]  (0,2.8) -- (0.6, 2.1) -- (1.5,1.7) -- (1.5, 1) -- (0,1);
		\fill [pattern=north east lines, pattern color=blue]  (0,2.8) -- (0.6, 2.1) -- (1.5,1.7) -- (1.5,3.4) -- (0, 3.4);
		\node[below] at (1.5,0) {{$\frac{n}{2}$}};
		\node[left] at (0,1) {{$1$}};
		\node[left] at (0,1.7) {{$p_{\mathrm{Fuj}}(n)$}};
		\draw[dashed, color=black]  (0, 1)--(3.6, 1);
		\draw[dashed, color=black]  (0, 1.7)--(1.5, 1.7);
		\draw[dashed, color=black]  (1.5, 0)--(1.5, 3.4);
		\draw[->] (5.8,0) -- (9.8,0) node[below] {$\gamma$};
		\draw[->] (6,-0.2) -- (6,3.4) node[left] {$p$};
		\node[left] at (6,-0.2) {{$0$}};
		\node[left] at (6,2.7) {{$\frac{n}{n-2}$}};
		\node[below] at (7,0) {{${\tilde{\gamma}}$}};
		\node[below] at (8,0) {{$\frac{n}{2}$}};
		\node[left] at (6,1) {{$1$}};
		\node[left] at (6,1.5) {{$p_{\mathrm{Fuj}}(n)$}};
		\node[left, color=black] at (9.64,1.7) {{$\longleftarrow$ $p=1+\frac{2\gamma}{n}$}};
		\node[right, color=black] at (3.6,2) {{ $p=p_{\mathrm{crit}}(n,\gamma)$$\longrightarrow$}};
		\draw[dashed, color=black]  (6, 1)--(9.6, 1);
		\draw[dashed, color=black]  (8, 0)--(8, 3.4);
		\draw[dashed, color=black]  (6, 2.7)--(8, 2.7);
		\draw[dashed, color=black]  (7, 0)--(7, 1.5);
		\draw[dashed, color=black] (6,1.5)--(7,1.5);
		\draw[color=black] plot[smooth, tension=.7] coordinates {(6,2.4) (6.3,1.9) (7,1.5)};
		\draw[dashed, color=black] plot[smooth, tension=.7] coordinates {(6.2,3.2) (6.4,2.3) (7,1.5)};
		\node[below] at (6.2,3.2) {{$\nearrow$}};
		\node[below] at (7.5,3.2) {{$\nwarrow$}};
		\draw[color=black] plot[smooth, tension=.7] coordinates { (7,1.5) (8,2)};
		\draw[dashed, color=black] plot[smooth, tension=.7] coordinates { (7,1.5) (7.7,3.3)};
		\draw[dashed, color=black] plot[smooth, tension=.7] coordinates { (7,1.5) (7.5,1.4) (8,1.3)};
		\fill [pattern=north west lines, pattern color=red]  (6,2.4) -- (6.3,1.9) -- (7,1.5) -- (7.5,1.4) --(8,1.3) -- (8,1) -- (6, 1);
		\fill [pattern=north east lines, pattern color=blue]  (6,2.4) -- (6.3,1.9) -- (7,1.5) -- (8,2) -- (8,2.7) -- (6,2.7);
		\node[left] at (3,-1) {{Cases $n=1,2$}};
		\node[left] at (9.2,-1) {{Cases $n=3,4,5,6$}};
	\end{tikzpicture}
	\caption{Description of the critical exponent in the $\gamma-p$ plane}
	\label{imgg}
\end{figure}

\begin{remark}\label{Rem_Higher_Dimension}
	In Theorem \ref{Thm_GESDS_Lower}, the global existence result holds for some low dimensions due to the technical assumption $1<p\leqslant\frac{n}{n-2s}$ for $n>2s$. In order to derive the global existence result for higher dimensions, we may study Sobolev solutions by considering initial data belonging to Sobolev spaces with suitable higher regularity or even large regularity. One may see \cite{Palmieri-Reissig-2018}. Nevertheless, in those situations, the applications of fractional chain rule or fractional powers in the treatment of nonlinear terms will bring new lower bound restrictions for the exponent $p$.
\end{remark}
We now turn our consideration to a blow-up result in the subcritical case.
\begin{theorem}\label{Thm_Blow_up}
	Let $\gamma\in(0,\frac{n}{2})$ for $n\geqslant 1$. Let the exponent $p$ fulfil $1<p<p_{\mathrm{crit}}(n,\gamma)$. Let us assume non-negative initial data $(u_0,u_1)\in \dot{H}^{-\gamma}\times\dot{H}^{-\gamma}$ such that
	\begin{align}\label{Special}
		u_0(x)+u_1(x)\geqslant \epsilon_1\langle x\rangle^{-n(\frac{1}{2}+\frac{\gamma}{n})}(\log(\mathrm{e}+|x|))^{-1},
	\end{align}
	where $\epsilon_1>0$ is a fixed constant. Then, there is no global (in time) weak solution to \eqref{Semilinear_Damped_Waves}. Furthermore, the lifespan $T_{\epsilon,\mathrm{w}}$ of local (in time) weak solutions to the Cauchy problem \eqref{Semilinear_Damped_Waves} fulfils
	\begin{align*}
		T_{\epsilon,\mathrm{w}}\leqslant C\epsilon^{-\frac{2}{2p'-2-\frac{n}{2}-\gamma}},
	\end{align*}
	where $C$ is a positive constant independent of $\epsilon$ and $p'$ is the conjugate exponent to $p$.
\end{theorem}
\begin{remark}
	The set of initial data, i.e. $(u_0,u_1)\in\dot{H}^{-\gamma}\times H^{-\gamma}$ with the condition \eqref{Special}, is non-empty, whose explanation will be shown in Section \ref{Section_Blow-up}.
\end{remark}
To end this subsection, we next analyze the critical exponent in the framework of Sobolev spaces of negative order. According to Theorems \ref{Thm_GESDS_Lower} and \ref{Thm_Blow_up}, we conclude the critical exponent for the semilinear damped wave equation \eqref{Semilinear_Damped_Waves} with initial data belonging additionally to $\dot{H}^{-\gamma}$ is \eqref{Critical_Exp}.
It provides a new viewpoint for the critical exponent of semilinear classical damped waves. Taking formally $\gamma=\frac{n}{2}$, the well-known Fujita exponent occurs such that $p_{\mathrm{crit}}(n,\frac{n}{2})=p_{\mathrm{Fuj}}(n)$. Or one may represent our critical exponent as $p_{\mathrm{crit}}(n,\gamma)=p_{\mathrm{Fuj}}(\frac{n}{2}+\gamma)$.
\begin{remark}
	We conjecture that by developing some suitable $\dot{H}^{-\gamma}_m-L^q$ estimates with $1\leqslant m\leqslant q\leqslant\infty$ for solutions to the corresponding linear classical damped wave equation with vanishing right-hand side, one may improve the global (in time) existence result when $\gamma>\tilde{\gamma}=\frac{\sqrt{n^2+16n}}{4}-\frac{n}{4}$ so that the restriction of $\gamma$ in the critical exponent \eqref{Critical_Exp} can be relaxed. Moreover, it remains open the question either to prove the global (in time) existence of small data Sobolev solutions or the blow-up of weak solutions in the critical case $p=p_{\mathrm{crit}}(n,\gamma)$.
\end{remark}
\noindent We take an example to show the critical exponent by considering Theorem \ref{Thm_GESDS_Lower} with $s=1$ and Theorem \ref{Thm_Blow_up}. Since the technical restriction from applications of fractional Gagliardo-Nirenberg inequality, we have to consider $1<p\leqslant \frac{n}{n-2}$ by Theorem \ref{Thm_GESDS_Lower} which leads to the empty set of admissible exponents $p$ if $n\geqslant 7$. Consequently, the ranges of exponent $p$ for global existence and blow-up are
\begin{itemize}
	\item when $n=1,2$: blow-up of weak solutions if $1<p<p_{\mathrm{crit}}(n,\gamma)$, and global (in time) existence of Sobolev solutions if $p>p_{\mathrm{crit}}(n,\gamma)$ for $0<\gamma<\frac{n}{2}$;
	\item when $n=3,4,5,6$: blow-up of weak solutions if $1<p<p_{\mathrm{crit}}(n,\gamma)$, and global (in time) existence of Sobolev solutions if $p_{\mathrm{crit}}(n,\gamma)<p\leqslant \frac{n}{n-2}$ for $0<\gamma\leqslant \tilde{\gamma}$, as well as $1+\frac{2\gamma}{n}\leqslant p\leqslant\frac{n}{n-2}$ for $\tilde{\gamma}<\gamma<\frac{n}{2}$.
\end{itemize}
\newpage
Precisely,  the analysis in this example can be described by Figure \ref{imgg}.

\noindent We remark that with an increasing  dimension  the curve $p=p_{\mathrm{crit}}(n,\gamma)$ and the segment $p=1+\frac{2\gamma}{n}$ will move following the direction of (dashed lines) arrows, respectively, in the second graph.

The statements of Theorems \ref{Thm_GESDS_Lower} and \ref{Thm_Blow_up} also show the critical index for regularity of initial data belonging additionally to $\dot{H}^{-\gamma}$ with $0<\gamma<\min\{ \tilde{\gamma},\frac{n}{2}\}$. To be specific, if $1+\frac{2}{n}<p<1+\frac{4}{n}$, the critical index is given by $\gamma_{\mathrm{crit}}(p,n):=\frac{2}{p-1}-\frac{n}{2}$. Namely, concerning initial data belonging additionally to $\dot{H}^{-\gamma}$, if $0<\gamma<\gamma_{\mathrm{crit}}(p,n)$, then local (in time) weak solutions in general blow up in finite time, and if $\gamma_{\mathrm{crit}}(p,n)<\gamma<\min\{ \tilde{\gamma},\frac{n}{2}\}$, then global (in time) small data Sobolev solutions exists uniquely.
\begin{remark}
	We expect our approach can be generalized to investigate the critical exponent for semilinear structurally damped evolution equations (even with some time-dependent coefficients $a(t)$, $b(t)$ in the ``effective'' case \cite{Wirth=2007,D'Abb-Lucente-Reissig=2013}), namely, for the Cauchy problem
	\begin{align*}
		\begin{cases}
			u_{tt}+a(t)(-\Delta)^{\sigma}u+b(t)(-\Delta)^{\delta} u_t=\begin{cases}
				|u|^p,\\
				|u_t|^p,
			\end{cases}&x\in\mb{R}^n, \ t>0,\\
			u(0,x)=u_0(x),\ u_t(0,x)=u_1(x),&x\in\mb{R}^n,
		\end{cases}
	\end{align*}
	with $p>1$, $\sigma\geqslant 1$, $\delta\in[0,\frac{\sigma}{2}]$, and the critical curve for its corresponding weakly coupled systems (see \cite{Pham-Kainane-Reissig=2015,D'Abbicco-Ebert=2017} and references therein), where the initial data belongs additionally to $\dot{H}^{-\gamma}$ for some $\gamma>0$.
\end{remark}

\subsection{Results and discussions for sharp lifespan estimates}
As we explained in Theorem \ref{Thm_Blow_up}, under the condition $1<p<p_{\mathrm{crit}}(n,\gamma)$, non-trivial local (in time) weak solution may blow up in finite time, which motivates us to provide more detailed information of the lifespan.

According to the derived result in Theorem \ref{Thm_Blow_up}, we have obtained upper bound estimates for the lifespan as follows:
\begin{align*}
	T_{\epsilon}\leqslant T_{\epsilon,\mathrm{w}}\leqslant C\epsilon^{-\frac{2}{2p'-2-\frac{n}{2}-\gamma}}.
\end{align*}
For this reason, the remaining part is to estimate $T_{\epsilon}$ from the below. Before studying this, let us define mild solutions to the Cauchy problem \eqref{Semilinear_Damped_Waves} on $[0,T)$ with $T>0$ for $u\in\ml{C}([0,T),H^1)$ as solutions of the operator equation
\begin{align}\label{mild_solution}
	u(t,x)=\epsilon E_0(t,x)\ast_{(x)}u_0(x)+\epsilon E_1(t,x)\ast_{(x)}u_1(x)+\int_0^tE_1(t-\tau,x)\ast_{(x)}|u(\tau,x)|^p\mathrm{d}\tau
\end{align}
for $t\in[0,T)$. In the above, $E_0=E_0(t,x)$ and $E_1=E_1(t,x)$ are the fundamental solutions to the linear Cauchy problem \eqref{Damped_Waves} with initial data $(v_0,v_1)=(\delta_0,0)$ and $(v_0,v_1)=(0,\delta_0)$, respectively. Here, $\delta_0$ is the Dirac distribution at $x=0$ with respect to the spatial variables.

Let us introduce by $T_{\epsilon,\mathrm{m}}$ the lifespan of a mild solution $u=u(t,x)$. Then, we have the next result.\bigskip

\begin{theorem}\label{Thm_Lower_Bound}
	Let $\gamma\in(0,\min\{2,\frac{n}{2}\})$ for $n\geqslant1$. Let the power exponent $p$ fulfil  $1<p<p_{\mathrm{crit}}(n,\gamma)$ such that
	\begin{align}\label{Restriction}
		1+\frac{2\gamma}{n}\leqslant p\leqslant\frac{n}{(n-2)_+}.
	\end{align}
	Let us assume $(u_0,u_1)\in\ml{A}_1$. Then, there exists a constant $\epsilon_2$ such that for any $\epsilon\in(0,\epsilon_2]$ the lifespan $T_{\epsilon,\mathrm{m}}$ of mild solutions $u=u(t,x)$ to the Cauchy problem \eqref{Semilinear_Damped_Waves} fulfils
	\begin{align}\label{Lower_Bound}
		T_{\epsilon,\mathrm{m}}\geqslant  C\epsilon^{-\frac{2}{2p'-2-\frac{n}{2}-\gamma}},
	\end{align}
	where $C$ is a positive constant independent of $\epsilon$, and depends on $p,n,\gamma$ as well as $\|(u_0,u_1)\|_{\ml{A}_1}$ only.
\end{theorem}
\begin{remark}
	Due to the facts that
	\begin{align*}
		1+\frac{2\gamma}{n}\leqslant\frac{n}{n-2}=1+\frac{2}{n-2}\ \ \mbox{and}\ \ 1+\frac{2\gamma}{n}<1+\frac{4}{n+2\gamma},
	\end{align*}
	for any fixed dimensions $n=n_0$ and the assumed admissible $\gamma$, we always can take a suitably small parameter $\gamma>0$ such that the set of $p$ in Theorem \ref{Thm_Lower_Bound} will be not empty.
\end{remark}
\begin{remark}
	Since $u$ in Theorem \ref{Thm_GESDS_Lower} is a mild solution to \eqref{Semilinear_Damped_Waves}, then by a density argument (see, for example, Proposition 3.1 in \cite{Ikeda-Wakasugi=2013}) this mild solution is also a weak solution to \eqref{Semilinear_Damped_Waves} with $T_{\epsilon,\mathrm{m}}\leqslant T_{\epsilon}$. Consequently, putting the statements of Theorems \ref{Thm_Blow_up} and \ref{Thm_Lower_Bound} together, we claim a sharp estimate for the lifespan $T_{\epsilon}$ when $1<p<p_{\mathrm{crit}}(n,\gamma)$, that is,
	\begin{align*}
		T_{\epsilon}\simeq C\epsilon^{-\frac{2}{2p'-2-\frac{n}{2}-\gamma}}=C\epsilon^{-\frac{2(p-1)}{2-(\frac{n}{2}+\gamma)(p-1)}}
	\end{align*}
	for $\gamma\in(0,\min\{2,\frac{n}{2}\})$. Furthermore, the lifespan estimate in the above coincides with the one for additionally $L^1$ regular data, i.e. with \eqref{Life_span}, if we formally take $\gamma=\frac{n}{2}$.
\end{remark}
\begin{remark}
	It is interesting to know the admissible ranges for $p$ in getting sharp lifespan estimates, which can be summarized as follows:
	\begin{itemize}
		\item when $n=1,2$, for  $1+\frac{2\gamma}{n}\leqslant p<p_{\mathrm{crit}}(n,\gamma)$ and $0<\gamma<\frac{n}{2}$, sharp lifespan estimates can be achieved;
		\item when $n=3,4$, for  $1+\frac{2\gamma}{n}\leqslant p<p_{\mathrm{crit}}(n,\gamma)$ and $0<\gamma\leqslant\frac{-n+\sqrt{n^2+16n}}{4}$, sharp lifespan estimates can be achieved;
		\item when $n=5,6,7$, for $1+\frac{2\gamma}{n}\leqslant p<p_{\mathrm{crit}}(n,\gamma)$ with $p\leqslant \frac{n}{n-2}$ and $0<\gamma\leqslant\frac{-n+\sqrt{n^2+16n}}{4}$, sharp lifespan estimates can be achieved;
		\item when $n\geqslant 8$, for $1+\frac{2\gamma}{n}\leqslant p\leqslant \frac{n}{n-2}$ and $0<\gamma\leqslant\frac{n}{n-2}$, sharp lifespan estimates can be achieved.
	\end{itemize}
\end{remark}
\begin{remark}
	The restriction \eqref{Restriction} of the range of $p$ for lower bound estimates \eqref{Lower_Bound}, again, comes from applications of the fractional Gagliardo-Nirenberg inequality. Namely for $n\geqslant 3$, we will pose an additional restriction $1<p\leqslant\frac{n}{n-2}$. It is also interesting to generalize those lower bound estimates without such upper bound for $p$ by constructing suitable time-spatial weighted Sobolev spaces in the derivation of lifespan estimates.
\end{remark}

\color{black}
\section{Linear damped wave equation with data belonging additionally to Sobolev spaces of negative order}\label{Section_Linear}
As a preparation to study the nonlinear Cauchy problem \eqref{Semilinear_Damped_Waves}, in this section, we consider the corresponding linear Cauchy problem with vanishing right-hand side, namely,
\begin{align}\label{Damped_Waves}
	\begin{cases}
		v_{tt}-\Delta v+v_t=0,&x\in\mb{R}^n,\ t>0,\\
		v(0,x)=v_0(x),\ v_t(0,x)=v_1(x),&x\in\mb{R}^n,
	\end{cases}
\end{align}
with initial data belonging additionally to Sobolev spaces of negative order. Particularly, we focus on decay properties and diffusion phenomenon of solutions, in which our proofs are slightly different from the classical theory. Our main tools are  asymptotic expansions and WKB analysis.

\subsection{Decay estimates of solutions}
First of all, let us apply the partial Fourier transform with respect to spatial variables to the Cauchy problem \eqref{Damped_Waves}. Then, we arrive at
\begin{align}\label{Damped_Waves_Fourier}
	\begin{cases}
		\hat{v}_{tt}+|\xi|^2 \hat{v}+\hat{v}_t=0,&\xi\in\mb{R}^n,\ t>0,\\
		\hat{v}(0,\xi)=\hat{v}_0(\xi),\ \hat{v}_t(0,\xi)=\hat{v}_1(\xi),&\xi\in\mb{R}^n,
	\end{cases}
\end{align}
whose characteristic equation is given by $\lambda^2+\lambda+|\xi|^2=0$. Thus, the eigenvalues are expressed by
\begin{align*}
	\lambda_{1,2}=\frac{1}{2}\big(-1\pm\sqrt{1-4|\xi|^2}\big).
\end{align*} By using asymptotic expansions of eigenvalues, we conclude the following relations:
\begin{itemize}
	\item $\lambda_1=-1+\ml{O}(|\xi|^2)$, $\lambda_2=-|\xi|^2+\ml{O}(|\xi|^4)$ for $|\xi|< \varepsilon\ll 1$;
	\item $\lambda_{1,2}=\pm i|\xi|-\frac{1}{2}+\ml{O}(|\xi|^{-1})$ for $|\xi|> N\gg 1$;
	\item $\Re\lambda_{1,2}>0$ for $\varepsilon\leqslant|\xi|\leqslant N$.
\end{itemize}
The solution to \eqref{Damped_Waves_Fourier} can be represented by
\begin{align*}
	\hat{v}(t,\xi)=\widehat{K}_0(t,|\xi|)\hat{v}_0(\xi)+\widehat{K}_1(t,|\xi|)\hat{v}_1(\xi),
\end{align*}
where the kernels in the Fourier space are
\begin{align*}
	\widehat{K}_0(t,|\xi|)&=\frac{\lambda_1\mathrm{e}^{\lambda_2t}-\lambda_2\mathrm{e}^{\lambda_1t}}{\lambda_1-\lambda_2}\\
	&=\begin{cases}
		\displaystyle{\frac{(-1+\ml{O}(|\xi|^2))\mathrm{e}^{(-|\xi|^2+\ml{O}(|\xi|^4))t}-(-|\xi|^2+\ml{O}(|\xi|^4))\mathrm{e}^{(-1+\ml{O}(|\xi|^2))t}}{-1+\ml{O}(|\xi|^2)}}&\mbox{for}\ \ |\xi|<\varepsilon,\\[1em]
		\displaystyle{\frac{(i|\xi|-\frac{1}{2}+\ml{O}(|\xi|^{-1}))\mathrm{e}^{(-i|\xi|-\frac{1}{2}+\ml{O}(|\xi|^{-1}))t}}{2i|\xi|+\ml{O}(1)}}&\\[1em]
		\qquad-\displaystyle{\frac{(-i|\xi|-\frac{1}{2}+\ml{O}(|\xi|^{-1}))\mathrm{e}^{(i|\xi|-\frac{1}{2}+\ml{O}(|\xi|^{-1}))t}}{2i|\xi|+\ml{O}(1)}}&\mbox{for}\ \ |\xi|> N,
	\end{cases}
\end{align*}
as well as
\begin{align*}
	\widehat{K}_1(t,|\xi|)=\frac{\mathrm{e}^{\lambda_1t}-\mathrm{e}^{\lambda_2t}}{\lambda_1-\lambda_2}=\begin{cases}
		\displaystyle{\frac{\mathrm{e}^{(-1+\ml{O}(|\xi|^2))t}-\mathrm{e}^{(-|\xi|^2+\ml{O}(|\xi|^4))t}}{-1+\ml{O}(|\xi|^2)}}&\mbox{for}\ \ |\xi| <\varepsilon,\\[1em]
		\displaystyle{\frac{\mathrm{e}^{(i|\xi|-\frac{1}{2}+\ml{O}(|\xi|^{-1}))t}-\mathrm{e}^{(-i|\xi|-\frac{1}{2}+\ml{O}(|\xi|^{-1}))t}}{2i|\xi|+\ml{O}(1)}}&\mbox{for}\ \ |\xi|> N.
	\end{cases}
\end{align*}
\begin{prop}\label{Prop_Decay_Linear}
	Let $(v_0,v_1)\in (H^s\cap\dot{H}^{-\gamma})\times(H^{s-1}\cap\dot{H}^{-\gamma})$ with $s\geqslant0$ and $s+\gamma\geqslant0$. Then, the solution $v=v(t,x)$ to the linear Cauchy problem \eqref{Damped_Waves} satisfies
	\begin{align}\label{Est_Hs}
		\|v(t,\cdot)\|_{\dot{H}^s}\lesssim (1+t)^{-\frac{s+\gamma}{2}}\big(\|v_0\|_{H^s\cap \dot{H}^{-\gamma}}+\|v_1\|_{H^{s-1}\cap \dot{H}^{-\gamma}}\big).
	\end{align}
\end{prop}
\begin{remark}
	Comparing with previous researches on wave equations with different damping mechanisms (see \cite{D'Abb-Reissig=2014,Ikehata=2014,Pham-Kainane-Reissig=2015} and references therein), we assume initial data belonging additionally to homogeneous Sobolev spaces with negative index $\dot{H}^{-\gamma}$ instead of additional $L^m$ regularity ($1\leqslant m<2$). Furthermore, the additional use of Sobolev spaces of negative order for initial data provides a decay rate $(1+t)^{-\frac{\gamma}{2}}$ for any $\gamma>0$.
\end{remark}
\begin{proof}
	From the previous asymptotic expansions of kernels, we can get some pointwise estimates in the Fourier space, namely,
	\begin{align*}
		|\widehat{K}_0(t,|\xi|)|&\lesssim\begin{cases}
			|\xi|^2\mathrm{e}^{-ct}+\mathrm{e}^{-c|\xi|^2t}&\mbox{for}\ \ |\xi|< \varepsilon\ll 1,\\
			\mathrm{e}^{-ct}&\mbox{for}\ \ \varepsilon\leqslant|\xi|\leqslant N,\\
			\mathrm{e}^{-ct}&\mbox{for}\ \ |\xi|> N\gg 1,
		\end{cases}\\
		|\widehat{K}_1(t,|\xi|)|&\lesssim\begin{cases}
			\mathrm{e}^{-ct}+\mathrm{e}^{-c|\xi|^2t}&\mbox{for}\ \ |\xi|< \varepsilon\ll 1,\\
			\mathrm{e}^{-ct}&\mbox{for}\ \ \varepsilon\leqslant|\xi|\leqslant N,\\
			|\xi|^{-1}\mathrm{e}^{-ct}&\mbox{for}\ \ |\xi|> N\gg 1,
		\end{cases}
	\end{align*}
	with suitable constants $c>0$. Concerning the small frequencies part, we apply
	\begin{align*}
		\|\,|\xi|^{s}\hat{f}(\xi)\hat{g}(\xi)\|_{L^2(|\xi|<\varepsilon)}\leqslant \|\,|\xi|^{s+\gamma}\hat{f}(\xi)\|_{L^{\infty}(|\xi|<\varepsilon)}\|\,|\xi|^{-\gamma}\hat{g}(\xi)\|_{L^2},
	\end{align*}
	associated with  Parseval's identity and
	\begin{align*}
		\|\,|\xi|^{s+\gamma}\mathrm{e}^{-c|\xi|^2t}\|_{L^{\infty}(|\xi|< \varepsilon)}\lesssim (1+t)^{-\frac{s+\gamma}{2}}\ \ \mbox{for}\ \ s+\gamma\geqslant0.
	\end{align*}
	Then, we have
	\begin{align}\label{Est_01}
		\|v(t,\cdot)\|_{\dot{H}^s}&=\big\|\,|\xi|^s(\widehat{K}_0(t,|\xi|)\hat{v}_0(\xi)+\widehat{K}_1(t,|\xi|)\hat{v}_1(\xi))\big\|_{L^2}\notag\\
		&\lesssim \big\|\,|\xi|^{s+\gamma}\widehat{K}_0(t,|\xi|)\big\|_{L^{\infty}(|\xi|<\varepsilon)}\|v_0\|_{\dot{H}^{-\gamma}}
		+\big\|\,|\xi|^{s+\gamma}\widehat{K}_1(t,|\xi|)\big\|_{L^{\infty}(|\xi|<\varepsilon)}\|v_1\|_{\dot{H}^{-\gamma}}\notag\\
		&\quad+\mathrm{e}^{-ct}(\|v_0\|_{L^2}+\|v_1\|_{H^{s-1}})+\|\widehat{K}_0(t,|\xi|)\|_{L^{\infty}(|\xi|>N)}\|v_0\|_{H^{s}}\notag\\
		&\quad+\|\widehat{K}_1(t,|\xi|)\|_{L^{\infty}(|\xi|>N)}\|v_1\|_{H^{s-1}}\notag\\
		&\lesssim (1+t)^{-\frac{s+\gamma}{2}}\big(\|v_0\|_{H^{-\gamma}}+\|v_1\|_{H^{-\gamma}}\big)+\mathrm{e}^{-ct}\big(\|v_0\|_{H^s}+\|v_1\|_{H^{s-1}}\big).
	\end{align}
	Our proof is complete.
\end{proof}
\subsection{Diffusion phenomenon}
It is well-known that the diffusion phenomenon bridges decay properties of solutions to the Cauchy problem for the classical damped wave equation \eqref{Damped_Waves} and solutions to the following Cauchy problem for the heat equation:
\begin{align}\label{Heats}
	\begin{cases}
		w_t-\Delta w=0,&x\in\mb{R}^n,\ t>0,\\
		w(0,x)=v_0(x)+v_1(x),&x\in\mb{R}^n.
	\end{cases}
\end{align}
Basing on $L^m$ initial data, a lot of papers derive a diffusion phenomenon in the sense that by measuring the difference of Sobolev solutions of \eqref{Damped_Waves} and \eqref{Heats} in suitable norms, an additional time decay rate appears. It shows large-time approximation with $L^m$ data. We refer interested readers to \cite{Hsiao-Liu=1992,Karch=2000,Nishihara=2003} and references therein.

In this part, we establish the diffusion phenomenon with initial data belonging additionally to Sobolev spaces of negative order. The solution to \eqref{Heats} in the Fourier space can be written by
\begin{align*}
	\hat{w}(t,\xi)=\mathrm{e}^{-|\xi|^2t}(\hat{v}_0(\xi)+\hat{v}_1(\xi)).
\end{align*}
Therefore, following the same procedure as in the proof of Proposition \ref{Prop_Decay_Linear}, it yields
\begin{align}\label{Est_H_Hs}
	\|w(t,\cdot)\|_{\dot{H}^s}\lesssim (1+t)^{-\frac{s+\gamma}{2}}\big(\|v_0\|_{ H^s\cap \dot{H}^{-\gamma}}+\|v_1\|_{H^s\cap \dot{H}^{-\gamma}}\big).
\end{align}
Comparing \eqref{Est_Hs} and \eqref{Est_H_Hs}, the solutions for the damped wave equation \eqref{Damped_Waves} and heat equation \eqref{Heats} fulfil the same decay estimates with slightly different regularity of data. It motivates us to study decay properties for the difference $v(t,\cdot)-w(t,\cdot)$ in $\dot{H}^s$.
\begin{prop}\label{Prop_Diffusion_Phenom}
	Let $(v_0,v_1)\in (H^s\cap\dot{H}^{-\gamma})\times(H^{s}\cap\dot{H}^{-\gamma})$ with $s\geqslant0$ and $s+\gamma+2\geqslant 0$. Then, the difference of solutions to the linear Cauchy problems \eqref{Damped_Waves} and \eqref{Heats} satisfies
	\begin{align}\label{Est_Diff_Hs}
		\|v(t,\cdot)-w(t,\cdot)\|_{\dot{H}^s}\lesssim (1+t)^{-\frac{s+\gamma}{2}-1}\big(\|v_0\|_{H^s\cap \dot{H}^{-\gamma}}+\|v_1\|_{H^{s}\cap \dot{H}^{-\gamma}}\big).
	\end{align}
\end{prop}
\begin{remark}
	Concerning the decay rate in \eqref{Est_Diff_Hs}, we observe the improvement $(1+t)^{-1}$ when we subtract the solution to the heat equation \eqref{Heats}. It concludes that diffusion phenomenon is also valid in the framework of $\dot{H}^{-\gamma}$ data. As a result, we may expect the sharpness of the derived estimate in Proposition \ref{Prop_Decay_Linear}.
\end{remark}
\begin{proof}
	To begin with, we define the cut-off function $\chi_{\intt}(\xi)\in\mathcal{C}^{\infty}$ with its support in $\{\xi\in\mb{R}^n:|\xi|<\varepsilon\ll 1\}$. According to the representation of solutions, we directly compute
	\begin{align*}
		|\hat{v}(t,\xi)-\hat{w}(t,\xi)|&\lesssim\chi_{\intt}(\xi)\Big(\Big|
		\frac{(-1+\ml{O}(|\xi|^2))\mathrm{e}^{(-|\xi|^2+\ml{O}(|\xi|^4))t}}{-1+\ml{O}(|\xi|^2)}-\mathrm{e}^{-|\xi|^2t}\Big|+\mathrm{e}^{-ct}\Big)|\hat{v}_0(\xi)|\\
		&\quad +\chi_{\intt}(\xi)\Big(\Big|-\frac{\mathrm{e}^{(-|\xi|^2+\ml{O}(|\xi|^4))t}}{-1+\ml{O}(|\xi|^2)}-\mathrm{e}^{-|\xi|^2t}\Big|+\mathrm{e}^{-ct}\Big)
		|\hat{v}_1(\xi)|\\
		&\quad+(1-\chi_{\intt}(\xi))\mathrm{e}^{-ct}|\hat{v}_0(\xi)|+(1-\chi_{\intt}(\xi))(1+|\xi|^{-1})\mathrm{e}^{-ct}|\hat{v}_1(\xi)|.
	\end{align*}
	From the fact
	\begin{align*}
		\chi_{\intt}(\xi)\big(\mathrm{e}^{(-|\xi|^2+\ml{O}(|\xi|^4))t}-\mathrm{e}^{-|\xi|^2t}\big)&=\chi_{\intt}(\xi)\mathrm{e}^{-|\xi|^2t}\ml{O}(|\xi|^4)t\int_0^1\mathrm{e}^{\ml{O}(|\xi|^4)ts}\mathrm{d}s\\
		&\lesssim\chi_{\intt}(\xi)|\xi|^2\mathrm{e}^{-c|\xi|^2t},
	\end{align*}
	we can deduce that the modulus of the difference of solutions to \eqref{Damped_Waves} and \eqref{Heats} satisfies
	\begin{align*}
		|\hat{v}(t,\xi)-\hat{w}(t,\xi)|&\lesssim\chi_{\intt}(\xi)\big(|\xi|^2\mathrm{e}^{-c|\xi|^2t}+\mathrm{e}^{-ct}\big)(|\hat{v}_0(\xi)|+|\hat{v}_1(\xi)|)\\
		&\quad+(1-\chi_{\intt}(\xi))\mathrm{e}^{-ct}(|\hat{v}_0(\xi)|+|\hat{v}_1(\xi)|).
	\end{align*}
	Following a similar procedure as those in the proof to Proposition \ref{Prop_Decay_Linear}, we complete the proof, in which the additional coefficient $|\xi|^2$ of the multiplier $\mathrm{e}^{-c|\xi|^2t}$ provides a faster decay rate.
\end{proof}
\section{Global (in time) well-posedness for semilinear Cauchy problems}\label{Section_GESDS}
\subsection{Philosophy of our approach}
For $T>0$, we introduce the evolution spaces
\begin{align*}
	X_s(T):=\ml{C}([0,T],H^s)\ \ \mbox{with}\ \ s\in(0,1],	
\end{align*}
carrying its norm
\begin{align*}
	\|u\|_{X_s(T)}:=\sup\limits_{t\in[0,T]}\big((1+t)^{\frac{\gamma}{2}}\|u(t,\cdot)\|_{L^2}+(1+t)^{\frac{s+\gamma}{2}}\|u(t,\cdot)\|_{\dot{H}^s}\big)
\end{align*}
with $\gamma>0$. The time-weighted Sobolev norm is strongly motivated by estimates for solutions to the linear Cauchy problem \eqref{Damped_Waves}, precisely, Proposition \ref{Prop_Decay_Linear} with $s=0$ and $s>0$, respectively. Then, let us define the operator $N$ such that
\begin{align*}
	N:\ u(t,x)\in X_s(T)\to Nu(t,x):=u^{\lin} (t,x)+u^{\non}(t,x),
\end{align*}
where $u^{\lin}(t,x) \equiv v(t,x)$ is the solution to the corresponding linear Cauchy problem \eqref{Damped_Waves} with the size $\epsilon$ for initial data, and $u^{\non}(t,x)$ is defined by
\begin{align*}
	u^{\non}(t,x):=\int_0^tE_1(t-\sigma,x)\ast_{(x)}|u(\sigma,x)|^p\mathrm{d}\sigma,
\end{align*}
where $E_1=E_1(t,x)$ is the fundamental solution to the linear Cauchy problem \eqref{Damped_Waves} with initial data $v_0(x)=0$ and $v_1(x)=\delta_0$. Here, $\delta_0$ is the Dirac distribution at $x=0$ with respect to the spatial variables.

In the subsequent part, we will prove global (in time) existence and uniqueness of small data Sobolev solutions of low regularity to the semilinear damped wave equation \eqref{Semilinear_Damped_Waves} by proving a unique fixed point of the operator $N$ that means $Nu\in X_s(T)$ for all positive $T$. To be specific, we will prove the following two crucial inequalities:
\begin{align}
	\|Nu\|_{X_s(T)}&\lesssim\epsilon\|(u_0,u_1)\|_{\ml{A}_{s}}+\|u\|_{X_s(T)}^p,\label{Cruc-01}\\
	\|Nu-N\bar{u}\|_{X_s(T)}&\lesssim\|u-\bar{u}\|_{X_s(T)}\big(\|u\|_{X_s(T)}^{p-1}+\|\bar{u}\|_{X_s(T)}^{p-1}\big),\label{Cruc-02}
\end{align}
respectively, with initial data space $\ml{A}_s:=(H^s\cap \dot{H}^{-\gamma})\times(L^2\cap \dot{H}^{-\gamma})$. In our aim inequality \eqref{Cruc-02}, $u$ and $\bar{u}$ are two solutions to the semilinear damped wave equation \eqref{Damped_Waves}. Providing that we take $\|(u_0,u_1)\|_{\ml{A}_{s} }\leqslant C$ and a small parameter $0<\epsilon\ll 1$, then we combine \eqref{Cruc-01} with \eqref{Cruc-02} to claim that there exists a global (in time) small data Sobolev solution $u^{*}=u^{*}(t,x)\in X_s(T)$ for all positive $T$ by using Banach's fixed point theorem.

To end this subsection, we recall two useful inequalities that will be applied later.
\begin{prop}\label{fractionalgagliardonirenbergineq} (Fractional Gagliardo-Nirenberg inequality, \cite{Hajaiej-Molinet-Ozawa-Wang-2011})
	Let $p,p_0,p_1\in(1,\infty)$ and $\kappa\in[0,s)$ with $s\in(0,\infty)$. Then, for all $f\in L^{p_0}\cap \dot{H}^{s}_{p_1}$ the following inequality holds:
	\begin{equation*}
		\|f\|_{\dot{H}^{\kappa}_{p}}\lesssim\|f\|_{L^{p_0}}^{1-\beta}\|f\|^{\beta}_{\dot{H}^{s}_{p_1}},
	\end{equation*}
	where $\beta=(\frac{1}{p_0}-\frac{1}{p}+\frac{\kappa}{n})\big/(\frac{1}{p_0}-\frac{1}{p_1}+\frac{s}{n})$ and $\beta\in[\frac{\kappa}{s},1]$.
\end{prop}
\begin{prop}\label{Hardy-Littlewood-Soblolev} (Hardy-Littlewood-Sobolev inequality, \cite{Lieb-1983})
	Let $0<\alpha<n$ and $1<m<q<\infty$ such that $\frac{1}{q}=\frac{1}{m}-\frac{\alpha}{n}$. Then, there exists a constant $C$ depending only on $m$ such that
	\begin{align*}
		\|\ml{I}_{\alpha}f\|_{L^q}\leqslant C\|f\|_{L^m},
	\end{align*}
	where $\ml{I}_{\alpha}f:=(-\Delta)^{-\frac{\alpha}{2}}f$ denotes the Riesz potential on $\mb{R}^n$.
\end{prop}
\subsection{Proof of Theorem \ref{Thm_GESDS_Lower}: Sobolev solutions of low regularity}
Before starting our proof, let us recall in Proposition \ref{Prop_Decay_Linear} two crucial estimates for Sobolev solutions to the linear Cauchy problem for $s\in[0,1]$ and $\gamma>0$ as follows:
\begin{align}
	\|v(t,\cdot)\|_{\dot{H}^s}&\lesssim (1+t)^{-\frac{s+\gamma}{2}}(\|v_0\|_{H^s\cap\dot{H}^{-\gamma}}+\|v_1\|_{L^2\cap\dot{H}^{-\gamma}}),\label{Eq_01}\\
	\|v(t,\cdot)\|_{\dot{H}^s}&\lesssim (1+t)^{-\frac{s}{2}}(\|v_0\|_{H^s}+\|v_1\|_{L^2}),\label{Eq_02}
\end{align}
where we used $L^2\hookrightarrow H^{s-1}$ for $s\leqslant 1$.

Since \eqref{Eq_01} and \eqref{Eq_02} with the size of initial data, we may claim $v\in X_s(T)$ such that
\begin{align*}
	\|u^{\lin}\|_{X_s(T)}\lesssim\epsilon\|(u_0,u_1)\|_{\ml{A}_s}.
\end{align*} Consequently, our first aim is to prove
\begin{align}\label{Crucial_03}
	\|u^{\non}\|_{X_s(T)}\lesssim\|u\|_{X_s(T)}^p
\end{align}
under some conditions for $p$. To do this, some estimates of $|u(\sigma,\cdot)|^p$ in $L^2$ and $\dot{H}^{-\gamma}$ are needed to be deduced. The fractional Gagliardo-Nirenberg inequality (see Proposition \ref{fractionalgagliardonirenbergineq}) implies
\begin{align*}
	\|\,|u(\sigma,\cdot)|^p\|_{L^2}\lesssim\|u(\sigma,\cdot)\|_{L^2}^{p(1-\beta_1)}\|u(\sigma,\cdot)\|_{\dot{H}^s}^{p\beta_1}\lesssim (1+\sigma)^{-(\frac{\gamma}{2}+\frac{n}{4})p+\frac{n}{4}}\|u\|^p_{X_s(T)}
\end{align*}
for $\sigma\in[0,T]$ with $\beta_1=\frac{n}{2s}(1-\frac{1}{p})\in[0,1]$. Therefore, we have to restrict $1<p\leqslant \frac{n}{n-2s}$ if $n>2s$. For another, we may use Hardy-Littlewood-Sobolev inequality (see Proposition \ref{Hardy-Littlewood-Soblolev}) so that
\begin{align*}
	\|\,|u(\sigma,\cdot)|^p\|_{\dot H^{-\gamma}}=\|\ml{I}_{\gamma}|u(\sigma,\cdot)|^p\|_{L^2}\lesssim\|\,|u(\sigma,\cdot)|^p\|_{L^m}=\|u(\sigma,\cdot)\|^p_{L^{mp}},
\end{align*}
with $\frac{1}{m}-\frac{1}{2}=\frac{\gamma}{n}$ carrying $\gamma\in(0,n)$ and $m\in(1,2)$. According to $m\in(1,2)$, we should restrict $\gamma\in(0,\frac{n}{2})$. Again, the application of the fractional Gagliardo-Nirenberg inequality shows
\begin{align*}
	\|\,|u(\sigma,\cdot)|^p\|_{\dot H^{-\gamma}}\lesssim\|u(\sigma,\cdot)\|_{L^2}^{p(1-\beta_2)}\|u(\sigma,\cdot)\|_{\dot{H}^s}^{p\beta_2}\lesssim(1+\sigma)^{-(\frac{\gamma}{2}+\frac{n}{4})p+\frac{\gamma}{2}+\frac{n}{4}}\|u\|^p_{X_s(T)}
\end{align*}
for $\sigma\in[0,T]$ with $\beta_2=\frac{n}{s}(\frac{1}{2}-\frac{1}{mp})\in[0,1]$. From the last condition of $\beta_2$, one finds
\begin{align*}
	p\geqslant \frac{2}{m}=1+\frac{2\gamma}{n},
\end{align*}
moreover
\begin{align}\label{Condi_01}
	p\leqslant\frac{2n}{m(n-2s)}\ \ \mbox{if}\ \ n>2s.
\end{align}
Since $m\in(1,2)$, the condition \eqref{Condi_01} can be guaranteed by $1<p\leqslant \frac{n}{n-2s}$. In conclusion, by taking
\begin{align*}
	1+\frac{2\gamma}{n}\leqslant p\begin{cases}
		<\infty&\mbox{if}\ \ n\leqslant 2s,\\[0.5em]
		\displaystyle{\leqslant\frac{n}{n-2s}}&\mbox{if}\ \ n> 2s,
	\end{cases}
\end{align*}
we assert for $\sigma \in [0,T]$ that
\begin{align}\label{Eq_03}
	\|\,|u(\sigma,\cdot)|^p\|_{L^2\cap \dot H^{-\gamma}}\lesssim (1+\sigma)^{-(\frac{\gamma}{2}+\frac{n}{4})p+\frac{\gamma}{2}+\frac{n}{4}}\|u\|^p_{X_s(T)}.
\end{align}
To estimate  the solution itself in $L^2$, we use the  $(L^2\cap \dot{H}^{-\gamma})-L^2$ estimate \eqref{Eq_01} in $[0,t]$ as well as the treatment of the nonlinear term \eqref{Eq_03}. Hence, we obtain
\begin{align*}
	\|u^{\non}(t,\cdot)\|_{L^2}&\lesssim\int_0^t(1+t-\sigma)^{-\frac{\gamma}{2}}(1+\sigma)^{-(\frac{\gamma}{2}
		+\frac{n}{4})p+\frac{\gamma}{2}+\frac{n}{4}}\,\mathrm{d}\sigma\|u\|^p_{X_s(T)}\\
	&\lesssim (1+t)^{-\frac{\gamma}{2}}\int_0^{\frac{t}{2}}(1+\sigma)^{-(\frac{\gamma}{2}+\frac{n}{4})p+\frac{\gamma}{2}+\frac{n}{4}}\,\mathrm{d}\sigma\|u\|^p_{X_s(T)}\\
	&\quad+(1+t)^{-(\frac{\gamma}{2}+\frac{n}{4})p+\frac{\gamma}{2}+\frac{n}{4}}\int_{\frac{t}{2}}^t(1+t-\sigma)^{-\frac{\gamma}{2}}\,\mathrm{d}\sigma \|u\|^p_{X_s(T)}.
\end{align*}
For the first integral, since the condition $p>1+\frac{4}{n+2\gamma}$, we can get the uniform integrability over $[0,\frac{t}{2}]$. Let us compute the second part precisely. So, we get
\begin{align*}
	\int_{\frac{t}{2}}^t(1+t-\sigma)^{-\frac{\gamma}{2}}\,\mathrm{d}\sigma\lesssim
	\begin{cases}
		(1+t)^{1-\frac{\gamma}{2}}&\mbox{if}\ \ \gamma<2,\\
		\ln(\mathrm{e}+t)&\mbox{if}\ \ \gamma=2,\\
		1&\mbox{if}\ \ \gamma>2.
	\end{cases}
\end{align*}
By considering $p>1+\frac{4}{n+2\gamma}$ if $\gamma\leqslant 2$ as well as $p>1+\frac{2\gamma}{n+2\gamma}$ if $\gamma>2$, we can get
\begin{align*}
	(1+t)^{-(\frac{\gamma}{2}+\frac{n}{4})p+\frac{\gamma}{2}+\frac{n}{4}}\int_{\frac{t}{2}}^t(1+t-\sigma)^{-\frac{\gamma}{2}}\,\mathrm{d}\sigma\lesssim (1+t)^{-\frac{\gamma}{2}}.
\end{align*}
Summarizing the last estimates, it holds
\begin{align*}
	(1+t)^{\frac{\gamma}{2}}	\|u^{\non}(t,\cdot)\|_{L^2}\lesssim \|u\|^p_{X_s(T)}.
\end{align*}
To estimate the solution itself in $\dot{H}^s$, we use the $(L^2\cap \dot{H}^{-\gamma})-\dot{H}^s$ estimate \eqref{Eq_01} in $[0,\frac{t}{2}]$ and the $L^2-\dot{H}^s$ estimate  \eqref{Eq_02} in $[\frac{t}{2},t]$. It means
\begin{align}\label{New_Hs}
	\|u^{\non}(t,\cdot)\|_{\dot{H}^s}&\lesssim(1+t)^{-\frac{s+\gamma}{2}}\int_0^{\frac{t}{2}}(1+\sigma)^{-(\frac{\gamma}{2}
		+\frac{n}{4})p+\frac{\gamma}{2}+\frac{n}{4}}\,\mathrm{d}\sigma\|u\|^p_{X_s(T)}\notag\\
	&\quad+(1+t)^{-(\frac{\gamma}{2}+\frac{n}{4})p+\frac{n}{4}}\int_{\frac{t}{2}}^t(1+t-\sigma)^{-\frac{s}{2}}\,\mathrm{d}\sigma\|u\|^p_{X_s(T)}.
\end{align}
Since $s\in(0,1]$ and $p>1+\frac{4}{n+2\gamma}$, we may claim
\begin{align*}
	(1+t)^{\frac{s+\gamma}{2}}\|u^{\non}(t,\cdot)\|_{\dot{H}^s}\lesssim \|u\|^p_{X_s(T)}.
\end{align*}
All in all, the desired estimate \eqref{Crucial_03} is proved.
\begin{remark}
	Here, we assume $p>1+\frac{4}{n+2\gamma}$ if $\gamma\leqslant 2$, and $p>1+\frac{2\gamma}{n+2\gamma}$ if $\gamma>2$ from the integrability and decay estimates of solution, moreover we need $p\geqslant 1+\frac{2\gamma}{n}$ from the application of the Gagliardo-Nirenberg inequality. However, when $\gamma>2$, we found
	\begin{align*}
		\max\Big\{1+\frac{2\gamma}{n+2\gamma}, 1+\frac{2\gamma}{n} \Big\}=1+\frac{2\gamma}{n}.
	\end{align*}
	So, the condition for the exponent $p$ is reduced to \eqref{Condition_p}.
\end{remark}
In order to prove \eqref{Cruc-02}, we notice that
\begin{align*}
	\|Nu-N\bar{u}\|_{X_s(T)}=\Big\| \int_0^tE_1(t-\sigma,\cdot)\ast_{(x)}\big(|u(\sigma,\cdot)|^p-|\bar{u}(\sigma,\cdot)|^p\big)\,\mathrm{d}\sigma \Big\|_{X_s(T)}.
\end{align*}
Thanks to  Hardy-Littlewood-Sobolev inequality, we may estimate
\begin{align*}
	\|\,|u(\sigma,\cdot)|^p-|\bar{u}(\sigma,\cdot)|^p\|_{\dot{H}^{-\gamma}}\lesssim\|\,|u(\sigma,\cdot)|^p-|\bar{u}(\sigma,\cdot)|^p\|_{L^m}
\end{align*}
with $\frac{1}{m}-\frac{1}{2}=\frac{\gamma}{n}$ equipping $\gamma\in(0,\frac{n}{2})$ and $m\in(1,2)$. The application of H\"older's inequality shows
\begin{align*}
	\|\,|u(\sigma,\cdot)|^p-|\bar{u}(\sigma,\cdot)|^p\|_{L^m}\lesssim\|u(\sigma,\cdot)-\bar{u}(\sigma,\cdot)\|_{L^{mp}}\big(\|u(\sigma,\cdot)\|^{p-1}_{L^{mp}}+\|\bar{u}(\sigma,\cdot)\|^{p-1}_{L^{mp}}\big).
\end{align*}
Finally, employing the fractional Gagliardo-Nirenberg inequality to estimate three terms on the right-hand side of the previous inequality, we can complete our desired estimate \eqref{Cruc-02}. Our proof is finished.

\section{Proof of Theorem \ref{Thm_Blow_up}: Blow-up of weak solutions}\label{Section_Blow-up}
We set $\varphi(x):=\langle x\rangle^{-n}$ for $n\geqslant 1$. Moreover, the test function $\eta=\eta(t)\in\ml{C}^{\infty}([0,\infty))$ with $\mathrm{supp}\,\eta(t)\subset[0,1]$ such that $\eta(t)=1$ if $t\in[0,\frac{1}{2}]$ and $\eta(t)=0$ if $t\in[1,\infty)$. For $R>0$, we define $\varphi_R(x):=\langle R^{-1}x\rangle^{-n}$ and $\eta_R(t):=\eta(R^{-2}t)$.

Assume by contradiction that $u=u(t,x)$ is a global (in time) weak solution to the Cauchy problem \eqref{Semilinear_Damped_Waves}. So, it holds that
\begin{align*}
	&\epsilon\int_{\mb{R}^n}(u_0(x)+u_1(x))\varphi_R(x)\,\mathrm{d}x+\int_0^{\infty}\int_{\mb{R}^n}|u(t,x)|^p\varphi_R(x)\eta_R(t)\,\mathrm{d}x\mathrm{d}t\\
	&\qquad=\int_0^t\int_{\mb{R}^n}u(t,x)(\partial_t^2-\Delta-\partial_t)(\varphi_R(x)\eta_R(t))\,\mathrm{d}x\mathrm{d}t.
\end{align*}
Let us introduce
\begin{align*}
	I_R:=\int_0^{\infty}\int_{\mb{R}^n}|u(t,x)|^p\varphi_R(x)\eta_R(t)\,\mathrm{d}x\mathrm{d}t.
\end{align*}
Because
\begin{align*}
	\int_{\mb{R}^n}(u_0(x)+u_1(x))^{\frac{2n}{n+2\gamma}}\,\mathrm{d}x&\geqslant \epsilon_1^{\frac{2n}{n+2\gamma}}\int_{\mb{R}^n}\langle x\rangle^{-n}(\log(\mathrm{e}+|x|))^{-\frac{2n}{n+2\gamma}}\,\mathrm{d}x\\
	&\geqslant C\int_0^{\infty}\langle r\rangle ^{-1}(\log(\mathrm{e}+r))^{-\frac{2n}{n+2\gamma}}\,\mathrm{d}r,
\end{align*}
as well as
\begin{align*}
	\|u_0+u_1\|_{L^{\frac{2n}{n+2\gamma}}}^{\frac{2n}{n+2\gamma}}\leqslant C\Big(\|u_0\|_{L^{\frac{2n}{n+2\gamma}}}^{\frac{2n}{n+2\gamma}}+\|u_1\|_{L^{\frac{2n}{n+2\gamma}}}^{\frac{2n}{n+2\gamma}}\Big)<\infty,
\end{align*}
we can assert if $(u_0,u_1)\in \big(L^{\frac{2n}{n+2\gamma}}\times L^{\frac{2n}{n+2\gamma}}\big)$ and \eqref{Special} holds, then
\begin{align*}
	C\int_0^{\infty}\langle r\rangle ^{-1}(\log(\mathrm{e}+r))^{-\frac{2n}{n+2\gamma}}\mathrm{d}r<\infty.
\end{align*}
Clearly, the last inequality is valid.
For this reason, the set
\begin{align*}
	\ml{D}_{n,\gamma}:=\Big\{(u_0,u_1)\, :\, u_0(x)+u_1(x)\geqslant \epsilon_1\langle x\rangle^{-n(\frac{1}{2}+\frac{\gamma}{n})}(\log(\mathrm{e}+|x|))^{-1}   \Big\}
\end{align*}
fulfils $(u_0,u_1)\in\ml{D}_{n,\gamma}\cap \big(L^{\frac{2n}{n+2\gamma}}\times L^{\frac{2n}{n+2\gamma}}\big)\neq\emptyset$ for $2n/(n+2\gamma)>1$. For example, if
\begin{align*}
	u_0(x)=u_1(x)=\epsilon_1\langle x\rangle^{-n(\frac{1}{2}+\frac{\gamma}{n})}(\log(\mathrm{e}+|x|))^{-1},
\end{align*} we can prove $u_0,u_1\in L^{\frac{2n}{n+2\gamma}}$. According to Hardy-Littlewood-Sobolev inequality, one gets $L^{\frac{2n}{n+2\gamma}}\subset\dot{H}^{-\gamma}$ since $\frac{n+2\gamma}{2n}-\frac{1}{2}=\frac{\gamma}{n}$ with $\gamma\in(0,\frac{n}{2})$. All in all, we may claim
$$(u_0,u_1)\in\ml{D}_{n,\gamma}\cap (\dot{H}^{-\gamma}\times \dot{H}^{-\gamma})\neq\emptyset$$ for $\gamma\in(0,\frac{n}{2})$. We now may apply the assumption \eqref{Special} to get
\begin{align*}
	\int_{\mb{R}^n}(u_0(x)+u_1(x))\varphi_R(x)\,\mathrm{d}x&\geqslant C\int_{|x|\leqslant R}(u_0(x)+u_1(x))\,\mathrm{d}x\\
	&\geqslant C\int_{|x|\leqslant R}\langle x\rangle ^{-n(\frac{1}{2}+\frac{\gamma}{n})}(\log(\mathrm{e}+|x|))^{-1}\,\mathrm{d}x\\
	&\geqslant C R^{\frac{n}{2}-\gamma}(\log R)^{-1}
\end{align*}
for $R\gg1$. On the other hand, taking into consideration some properties of the test functions, and applying Young's inequality, we obtain
\begin{align*}
	\int_0^{\infty}\int_{\mb{R}^n}\big|u(t,x)(\partial_t^2-\Delta-\partial_t)(\varphi_R(x)\eta_R(t))\big|\,\mathrm{d}x\mathrm{d}t\leqslant CI_R^{\frac{1}{p}}R^{\frac{n+2}{p'}-2}\leqslant\frac{I_R}{p}+\frac{C}{p'}R^{n+2-2p'},
\end{align*}
where we used the fact that
\begin{align*}
	&\int_0^{\infty}\int_{\mb{R}^n}|u(t,x)\partial_t\eta_R(t)\varphi_R(x)|\,\mathrm{d}x\mathrm{d}t\\
	&\qquad\leqslant C R^{-2}I_R^{\frac{1}{p}}\Big(\int_0^{\infty}\int_{\mb{R}^n}(\varphi_R(x)\eta_R(t))^{-\frac{p'}{p}}(\eta'(R^{-2}t)\varphi_R(x))^{p'}\,\mathrm{d}x\mathrm{d}t\Big)^{\frac{1}{p'}}\\
	&\qquad\leqslant CI_R^{\frac{1}{p}}R^{\frac{n+2}{p'}-2}.
\end{align*}
Now, we state the contradiction as follows:
\begin{align}\label{New_Blow_up_crucial}
	0\leqslant\frac{1}{p'}I_R\leqslant\frac{C}{p'}R^{n+2-2p'}-C\epsilon R^{\frac{n}{2}-\gamma}(\log R)^{-1}<0
\end{align}	
for $R\gg1$ providing that $R^{n+2-2p'}< CR^{\frac{n}{2}-\gamma}(\log R)^{-1}$, namely, we require $n+2-2p'<\frac{n}{2}-\gamma$ or $p<1+\frac{4}{n+2\gamma}$. This completes the proof of the blow-up result.
\begin{remark} To get upper bound estimates for the lifespan, we take $R\uparrow T^{\frac{1}{2}}_{\epsilon,\mathrm{w}}$ in \eqref{New_Blow_up_crucial}.
	So, a contradiction holds if
	\begin{align*}
		\frac{C}{p'}T_{\epsilon,\mathrm{w}}^{\frac{n}{2}+1-p'}<C\epsilon T_{\epsilon,\mathrm{w}}^{\frac{n}{4}-\frac{\gamma}{2}}(\log T_{\epsilon,\mathrm{w}})^{-1}.
	\end{align*}
	A direct computation yields
	\begin{align*}
		T_{\epsilon,\mathrm{w}}\leqslant C\epsilon^{-\frac{2}{2p'-2-\frac{n}{2}-\gamma}}
	\end{align*}
	to complete our desired lifespan estimate from above. \end{remark}
\section{Proof of Theorem \ref{Thm_Lower_Bound}: Sharp lifespan estimates}\label{Section_Lower_Bound}
In this section, we will use some notations from Section \ref{Section_GESDS}, especially, the evolution space $X_1(T)$ and data space $\ml{A}_1$. Nevertheless, to develop lower bound estimates for the lifespan, we rely on another nonlinear inequality rather than \eqref{Cruc-01}, because we consider now $1<p<p_{\mathrm{crit}}(n,\gamma)$.

First of all, thanks to the definition of mild solutions in \eqref{mild_solution}, we can estimate and represent local (in time) mild solutions directly. By the same procedure as proposed in Section \ref{Section_GESDS}, we can get
\begin{align*}
	(1+t)^{\frac{\gamma}{2}}\|u(t,\cdot)\|_{L^2}&\lesssim \epsilon\|(u_0,u_1)\|_{\ml{A}_1}+\int_0^{\frac{t}{2}}(1+\sigma)^{-(\frac{\gamma}{2}+\frac{n}{4})p+\frac{\gamma}{2}+\frac{n}{4}}\mathrm{d}\sigma\,\|u\|_{X_1(T)}^p\\
	&\quad+(1+t)^{-(\frac{\gamma}{2}+\frac{n}{4})p+\frac{\gamma}{2}+\frac{n}{4}+1}\|u\|_{X_1(T)}^p
\end{align*}
because of $\gamma\in(0,\min\{2,\frac{n}{2}\})$, where we restricted from the application of the Gagliardo-Nirenberg inequality that
\begin{align*}
	1+\frac{2\gamma}{n}\leqslant p\begin{cases}
		<\infty&\mbox{if}\ \ n\leqslant 2,\\[0.5em]
		\displaystyle{\leqslant\frac{n}{n-2}}&\mbox{if}\ \ n\geqslant 3.
	\end{cases}
\end{align*}
Due to the fact that
\begin{align*}
	1<p<p_{\mathrm{crit}}(n,\gamma)\ \ \mbox{implying}\ \ -\Big(\frac{\gamma}{2}+\frac{n}{4}\Big)p+\frac{\gamma}{2}+\frac{n}{4}>-1,
\end{align*}
one may gain
\begin{align*}
	(1+t)^{\frac{\gamma}{2}}\|u(t,\cdot)\|_{L^2}\lesssim \epsilon\|(u_0,u_1)\|_{\ml{A}_1}+ (1+t)^{-(\frac{\gamma}{2}+\frac{n}{4})p+\frac{\gamma}{2}+\frac{n}{4}+1}\|u\|_{X_1(T)}^p.
\end{align*}
Choosing $s=1$ in \eqref{New_Hs}, it also leads to
\begin{align*}
	(1+t)^{\frac{1+\gamma}{2}}\|u(t,\cdot)\|_{\dot{H}^1}\lesssim \epsilon\|(u_0,u_1)\|_{\ml{A}_1}+ (1+t)^{-(\frac{\gamma}{2}+\frac{n}{4})p+\frac{\gamma}{2}+\frac{n}{4}+1}\|u\|_{X_1(T)}^p.
\end{align*}
That is to say
\begin{align}\label{Crucial_4}
	\|u\|_{X_1(T)}\leqslant \epsilon C_0+C_1(1+t)^{\alpha_0(p,n,\gamma)}\|u\|_{X_1(T)}^p,
\end{align}
where $\alpha_0(p,n,\gamma):=-(\frac{\gamma}{2}+\frac{n}{4})p+\frac{\gamma}{2}+\frac{n}{4}+1$ belongs to $(0,1)$, and $C_0,C_1$ are two positive constants independent of $\epsilon,T$.

Let us now introduce
\begin{align*}
	T^*:=\sup\big\{ T\in[0,T_{\epsilon,\mathrm{m}})\ \ \mbox{such that}\ \ \ml{G}(T):=\|u\|_{X_1(T)}\leqslant M\epsilon \big\}
\end{align*}
with a sufficient large constant $M>0$ to be chosen later. Hence, from the fact $\ml{G}(T^*)\leqslant M\epsilon$ and \eqref{Crucial_4}, we arrive at
\begin{align}\label{New_Ineq_01}
	\ml{G}(T^*)\leqslant \epsilon\big(C_0+C_1(1+T^*)^{\alpha_0(p,n,\gamma)}M^p\epsilon^{p-1}\big)
	< \frac{1}{2}M\epsilon<M\epsilon
\end{align}
for large $M$ such that $4C_0< M$ as well as
\begin{align*}
	4C_1(1+T^*)^{\alpha_0(p,n,\gamma)}M^{p-1}\epsilon^{p-1}<1.
\end{align*}
Note that $\ml{G}=\ml{G}(T)$ is a continuous function for any $T\in(0,T_{\epsilon,\mathrm{m}})$. Nonetheless, \eqref{New_Ineq_01} shows that there exists a time $T_0\in(T^*,T_{\epsilon,\mathrm{m}})$ such that $\ml{G}(T_0)\leqslant M\epsilon$, which contradicts to the definition of $T^*$. In other words, we need to pose the next condition:
\begin{align*}
	4C_1(1+T^*)^{\alpha_0(p,n,\gamma)}M^{p-1}\epsilon^{p-1}\geqslant 1.
\end{align*}
It means that we can derive the blow-up time estimate to below
\begin{align*}
	T_{\epsilon,\mathrm{m}}\geqslant C\epsilon^{-\frac{p-1}{\alpha_0(p,n,\gamma)}}.
\end{align*}
In this way we achieve our aim of lower bound estimates of lifespan.

\color{black}
\section*{Acknowledgments}
This work was supported by the China Postdoctoral Science Foundation (Grant No. 2021T140450 and No. 2021M692084).

\end{document}